\documentclass[11pt]{article}
\usepackage{amssymb,amsmath,amsthm,fullpage}

 \newtheorem*{Theorem 11}{Theorem 11}
 \newtheorem*{Conjecture 7}{Conjecture 7}
 \newtheorem{thm}{Theorem}
 \newtheorem{cor}[thm]{Corollary}
 \newtheorem{lem}[thm]{Lemma}
 
 \newtheorem{conj}[thm]{Conjecture}
 \theoremstyle{definition}
 
 \theoremstyle{remark}

\begin{document}

\title{Some bounds on convex combinations of $\omega$ and $\chi$ for decompositions into many parts}

\author{ Landon Rabern \\
         Department of Mathematics, UC Santa Barbara \\
         landonr@math.ucsb.edu}

\date{\today}

\maketitle

\begin{abstract}
A \emph{$k$--decomposition} of the complete graph $K_n$ is a decomposition of $K_n$ into $k$ spanning subgraphs
$G_1,\ldots,G_k$. For a graph parameter $p$, let $p(k;K_n)$ denote the maximum of $\displaystyle \sum_{j=1}^{k}
p(G_j)$ over all $k$--decompositions of $K_n$.  It is known that $\chi(k;K_n) = \omega(k;K_n)$ for $k \leq 3$
and conjectured that this equality holds for all $k$.  In an attempt to get a handle on this, we study convex
combinations of $\omega$ and $\chi$; namely, the graph parameters $A_r(G) = (1-r) \omega(G) + r \chi(G)$ for $0
\leq r \leq 1$.  It is proven that $A_r(k;K_n) \leq n  + {k \choose 2}$ for small $r$.  In addition, we prove
some generalizations of a theorem of Kostochka, et al. \cite{kostochka}.
\end{abstract}

\section{Introduction}

A \emph{$k$--decomposition} of the complete graph $K_n$ is a decomposition of $K_n$ into $k$ spanning subgraphs
$G_1,\ldots,G_k$; that is, the $G_j$ have the same vertices as $K_n$ and each edge of $K_n$ belongs to precisely
one of the $G_j$.  For a graph parameter $p$ and a positive integer $k$, define

\[p(k;K_n) = \max \{ \displaystyle \sum_{j=1}^{k} p(G_j) \mid (G_1,\ldots,G_k) \text{ a $k$--decomposition of } K_n \}. \]

We say $(G_1,\ldots,G_k)$ is a \emph{$p$-optimal} $k$-decomposition of $K_n$ if $\displaystyle \sum_{j=1}^{k}
p(G_j) = p(k;K_n)$. We will be interested in parameters that are convex combinations of the clique number and
the chromatic number of a graph $G$. For $0 \leq r \leq 1$, define $A_r(G) = (1-r) \omega(G) + r \chi(G)$. We
would like to determine $A_r(k;K_n)$.  The following theorem of Kostochka, et al. does this for the case $r=0$.

\begin{thm}[Kostochka, et al. \cite{kostochka}]
If $k$ and $n$ are positive integers, then $\omega(k;K_n) \leq n + {k \choose 2}$.  If $n \geq {k \choose 2}$,
then $\omega(k;K_n) = n + {k \choose 2}$.
\end{thm}

Since $A_r(k;K_n) \leq (1-r) \omega(k;K_n) + r \chi(k;K_n)$, this theorem combined with the following result of
Watkinson gives the general upper bound

\begin{equation}\label{generalupperbound}
A_r(k;K_n) \leq n + (1-r){k \choose 2} + r\frac{k!}{2}.
\end{equation}

\begin{thm}[Watkinson \cite{watkinson}]
If $k$ and $n$ are positive integers, then $\chi(k;K_n) \leq n + \frac{k!}{2}$.
\end{thm}

From Theorem 1, we see that $A_r(k;K_n) \leq n + {k \choose 2}$ is the best possible bound. Equation
\eqref{generalupperbound} shows that this holds for $k \leq 3$. Also, this bound is an immediate consequence of
a conjecture made by Plesn\'{i}k.

\begin{conj}[Plesn\'{i}k \cite{plesnik}]
If $k$ and $n$ are positive integers, then $\chi(k;K_n) \leq n + {k \choose 2}$.
\end{conj}

Since $\omega \leq \chi$, if the conjectured bound on $A_r(k;K_n)$ holds for $r$, then it holds for all $0\leq s
\leq r$ as well. This suggests that it may be easier to look at small values of $r$ first.  Our next theorem
proves the optimal bound for small $r$.

\begin{Theorem 11}
Let $k$ and $n$ be positive integers and $0\leq r \leq \min \{ 1,3/k \}$. Then
\[A_r(k;K_n) \leq n  + {k \choose 2}.\]
\end{Theorem 11}

Along the way we prove some generalizations of Theorem 1.  A definition is useful here. For $0 \leq m \leq k$,
define

\[\chi_m(k;K_n) = \max \{ \displaystyle \sum_{j=1}^{m} \chi(G_j) + \sum_{j=m+1}^{k} \omega(G_j) \mid (G_1,\ldots,G_k) \text{ a $k$--decomposition of } K_n \}. \]

We say $(G_1,\ldots,G_k)$ is a \emph{$\chi_m$-optimal} $k$-decomposition of $K_n$ if $\displaystyle
\sum_{j=1}^{m} \chi(G_j) + \sum_{j=m+1}^{k} \omega(G_j) = \chi_m(k;K_n)$.  Note that $\chi_0(k;K_n) =
\omega(k;K_n)$ and $\chi_k(k;K_n) = \chi(k;K_n)$. \bigskip

We prove that the following holds for a given value of $m$ if and only if Conjecture 3 holds for $k=m$.
\begin{Conjecture 7}
Let $m$ and $n\geq 1$ be non-negative integers.  Then $\chi_m(k;K_n) \leq n + {k \choose 2}$ for all $k \geq m$.
\end{Conjecture 7}

In the last section, we prove similar results for decompositions of $K_n^r$ into $r$-uniform hypergraphs.

\section{Notation}
We quickly fix some terminology and notation. \newline A \emph{hypergraph} $G$ is a pair consisting of finite
set $V(G)$ together with a set $E(G)$ of subsets of $V(G)$ of size at least two.  The elements of $V(G)$ and
$E(G)$ are called \emph{vertices} and \emph{edges} respectively. If $|e| = r$ for all $e \in E(G)$, then $G$ is
\emph{$r$-uniform}.  A $2$-uniform hypergraph is a \emph{graph}.  The \emph{order} $|G|$ of $G$ is the number of
vertices in $G$.  The \emph{size} $s(G)$ of $G$ is the number of edges in $G$.  The \emph{degree} $d(v)$ of a
vertex $v \in V(G)$ is the number of edges of $G$ that contain $v$.  Vertices $v_1,\ldots,v_t$ are called
\emph{adjacent} in $G$ if $\{v_1,\ldots,v_t\} \in E(G)$.
\newline
Given two hypergraphs $G$ and $H$, we say that $H$ is a \emph{subhypergraph} of $G$ if $V(H) \subseteq V(G)$ and
$E(H) \subseteq E(G)$.\newline Given a hypergraph $G$ and $X \subseteq V(G)$, let $G[X]$ denote the hypergraph
with vertex set $X$ and edge set $\{ e \in E(G) \mid e \subseteq X \}$.  This is called the subhypergraph of $G$
\emph{induced} by $X$.  Let $G - X$ denote $G[V(G) \smallsetminus X]$.  For $e \subseteq V(G)$, let $G+e$ and
$G-e$ denote the hypergraphs with vertex set $V(G)$ and edge sets $E(G) \cup \{e\}$ and $E(G) \smallsetminus
\{e\}$ respectively.\newline Given an $r$-uniform hypergraph $G$, $X \subseteq V(G)$ is a \emph{clique} if
$E(G[X])$ contains every $r$-subset of $X$. The \emph{clique number} $\omega(G)$ is the maximum size of a clique
in $G$. If $\omega(G) = |G|$, then $G$ is called \emph{complete}.  Denote the $r$-uniform complete hypergraph on
$n$ vertices by $K_n^r$.  For the case of graphs($r=2$) we drop the superscript, writing $K_n$.\newline  For a
graph $G$, the \emph{chromatic number} $\chi(G)$ of $G$ is the least number of labels required to label the
vertices so that adjacent vertices receive distinct labels.  Note that if $\{X_1,\ldots,X_t\}$ is a partition of
$V(G)$, then $\chi(G) \leq \displaystyle \sum_{j=1}^{t} \chi(G[X_j])$.  Following \cite{kostochka} we call this
property \emph{subadditivity} of $\chi$.

\section{Convex combinations of $\omega$ and $\chi$}
Given a graph $G$, let $P(G)$ denote the induced subgraph of $G$ on the vertices of positive degree; that is,

\[P(G) = G[\{v \in V(G) \mid d(v) \geq 1 \}]. \]
\begin{lem}
Let $0 \leq m < k$ and $n$ a positive integer.  If $(G_1,\ldots,G_k)$ is a $\chi_m$-optimal $k$-decomposition of
$K_n$ with $s(G_1)$ maximal, then $P(G_j)$ is complete for $m < j \leq k$.
\end{lem}
\begin{proof}
Let $(G_1,\ldots,G_k)$ be a $\chi_m$-optimal $k$-decomposition of $K_n$ with $s(G_1)$ maximal.  Let $m < j \leq
k$.  Take $e \in E(G_j)$.  Then $(G_1 + e,\ldots,G_j - e,\ldots,G_k)$ is a $k$-decomposition of $K_n$ with
$s(G_1 + e) > s(G_1)$.  Hence $(G_1 + e,\ldots,G_j - e,\ldots,G_k)$ is not $\chi_m$-optimal, which implies that
$\omega(G_j - e) < \omega(G_j)$.  Whence every edge of $G_j$ is involved in every maximal clique and thus every
vertex of positive degree is involved in every maximal clique.  Hence $\omega(P(G_j)) = |P(G_j)|$, showing
$P(G_j)$ complete.
\end{proof}

\begin{thm}
Let $m \geq 1$. Assume $\chi(m;K_n) \leq n + f(m)$ for all $n \geq 1$. Then, for $k \geq m$,
\[\chi_m(k;K_n) \leq n + {k \choose 2} + f(m) - {m \choose 2}.\]
\end{thm}
\begin{proof}
Fix $k \geq m$.  Let $(G_1,\ldots,G_k)$ be a $\chi_m$-optimal $k$-decomposition of $K_n$ with $s(G_1)$ maximal.
Set $X = \displaystyle \bigcup_{j=m+1}^{k} V(P(G_j))$.  Then $(G_1 - X,\ldots, G_m - X)$ is an $m$-decomposition
of $K_{n-|X|}$ and hence
\begin{equation}\label{leftovers}
\displaystyle \sum_{j=1}^{m} \chi(G_j - X) \leq n - |X| + f(m).
\end{equation}

Fix $1 \leq j \leq m$.  By Lemma 4, $P(G_i)$ is complete for $i > m$.  Hence $P(G_j[X])$ and $P(G_i[X])$ have at
most one vertex in common for $i > m$.  Thus $|P(G_j[X])| \leq k-m$.  In particular, $\chi(G_j[X]) =
\chi(P(G_j[X])) \leq k-m$.  Combining this with \eqref{leftovers}, we have

\[ \displaystyle \sum_{j=1}^{m} \chi(G_j - X)  + \sum_{j=1}^{m} \chi(G_j[X]) \leq n - |X| + f(m) + m(k-m).\]

By subadditivity of $\chi$, this is

\begin{equation}\label{chichunk}
\displaystyle \sum_{j=1}^{m} \chi(G_j) \leq n - |X| + f(m) + m(k-m).
\end{equation}

Also, since $P(G_i)$ is complete for $i > m$,
\[\displaystyle \sum_{i=m+1}^{k} \omega(G_i) = \sum_{i=m+1}^{k} |P(G_i)| \leq |X| + {k-m \choose 2}.\]

Adding this to \eqref{chichunk} yields

\[ \chi_m(k;K_n) = \displaystyle \sum_{j=1}^{m} \chi(G_j) + \displaystyle \sum_{i=m+1}^{k} \omega(G_i) \leq n + {k-m \choose 2}
+ f(m) + m(k-m),\]

which is the desired inequality since ${k-m \choose 2}+ m(k-m) = {k \choose 2} - {m \choose 2}$.
\end{proof}

\begin{cor}
Let $m \geq 1$. Assume $\chi(m;K_n) \leq n + {m \choose 2}$ for all $n \geq 1$. Then, for $k \geq m$,
\[\chi_m(k;K_n) \leq n + {k \choose 2}.\]
\end{cor}

This shows that the following holds for a given value of $m$ if and only if Conjecture 3 holds for $k=m$.

\begin{conj}
Let $m$ and $n\geq 1$ be non-negative integers.  Then $\chi_m(k;K_n) \leq n + {k \choose 2}$ for all $k \geq m$.
\end{conj}

Since $\chi(1;K_n) \leq n$, we immediately have a generalization of Theorem 1.

\begin{cor}
If $k$ and $n$ are positive integers, then $\chi_1(k;K_n) \leq n + {k \choose 2}$. If $n \geq {k \choose 2}$
then $\chi_1(k;K_n) = n + {k \choose 2}$.
\end{cor}

With the help of Theorem 2, we get a stronger generalization.

\begin{cor}
If $k \geq 3$ and $n$ are positive integers, then $\chi_3(k;K_n) \leq n + {k \choose 2}$.  If $n \geq {k \choose
2}$ then $\chi_3(k;K_n) = n + {k \choose 2}$.
\end{cor}

We don't know if Conjecture 7 holds for any larger value of $m$.

\begin{cor}
Let $k$ and $n$ be positive integers with $n \geq {k \choose 2}$.  If $A$ is a graph appearing in an
$\omega$-optimal $k$-decomposition of $K_n$, then $\chi(A) = \omega(A)$.
\end{cor}
\begin{proof}
Let $(A,G_2,\ldots,G_k)$ be an $\omega$-optimal $k$-decomposition of $K_n$.  Then, by Theorem 1,

\[\omega(A) + \displaystyle \sum_{j=2}^{k} \omega(G_j) = n + {k \choose 2}.\]

Hence, by Corollary 8,

\[n + {k \choose 2} = \omega(A) + \displaystyle \sum_{j=2}^{k} \omega(G_j) \leq \chi(A) + \displaystyle \sum_{j=2}^{k} \omega(G_j) \leq n + {k \choose 2}.\]

Thus,

\[\omega(A) + \displaystyle \sum_{j=2}^{k} \omega(G_j) \leq \chi(A) + \displaystyle \sum_{j=2}^{k} \omega(G_j),\]

which gives $\chi(A) = \omega(A)$ as desired.

\end{proof}

\begin{thm}
Let $k$ and $n$ be positive integers and $0\leq r \leq \min \{1,3/k \}$. Then
\[A_r(k;K_n) \leq n  + {k \choose 2}.\]
\end{thm}
\begin{proof}
If $k \leq 3$, then $r=1$ and the assertion follows from Corollary 9.  Assume $k > 3$.  Let $(G_1,\ldots,G_k)$
be a $k$--decomposition of $K_n$. Since any rearrangement of $(G_1,\ldots,G_k)$ is also a $k$--decomposition of
$K_n$, Corollary 9 gives us the ${k \choose 3}$ permutations of the inequality

\[\chi(G_1)  +  \chi(G_2)  +  \chi(G_3) +  \omega(G_4) + \ldots + \omega(G_{k})  \leq n + {k \choose 2}.\]

Adding these together gives

\[ {k-1 \choose 3} \displaystyle \sum_{j=1}^{k} \omega(G_j) + {k-1 \choose 2} \displaystyle \sum_{j=1}^{k} \chi(G_j)\leq
{k \choose 3} \left(n  +  {k \choose 2}\right),\] which is

\[\frac{k-3}{k} \displaystyle \sum_{j=1}^{k} \omega(G_j) + \frac{3}{k} \displaystyle \sum_{j=1}^{k} \chi(G_j) \leq
n + {k \choose 2}. \]

Combining the sums yields

\[\displaystyle \sum_{j=1}^{k} A_r(G_j) \leq \sum_{j=1}^{k} A_{\frac{3}{k}}(G_j) = \sum_{j=1}^{k} \left(\frac{k-3}{k}\omega(G_j) + \frac{3}{k}\chi(G_j) \right) \leq n + {k \choose 2}. \]

\end{proof}

\section{Clique number of uniform hypergraphs}
A \emph{$k$--decomposition} of the complete $r$-uniform hypergraph $K_n^r$ is a decomposition of $K_n^r$ into
$k$ spanning subhypergraphs $G_1,\ldots,G_k$; that is, the $G_j$ have the same vertices as $K_n^r$ and each edge
of $K_n^r$ belongs to precisely one of the $G_j$.  Let

\[\omega(k;K_n^r) = \max \{ \displaystyle \sum_{j=1}^{k} \omega(G_j) \mid (G_1,\ldots,G_k) \text{ a $k$--decomposition of } K_n^r \}. \]

We say $(G_1,\ldots,G_k)$ is a \emph{$\omega$-optimal} $k$-decomposition of $K_n^r$ if $\displaystyle
\sum_{j=1}^{k} \omega(G_j) = \omega(k;K_n^r)$.

Given an $r$-uniform hypergraph $G$, let $P(G)$ denote the induced subhypergraph of $G$ on the vertices of
positive degree; that is,

\[P(G) = G[\{v \in V(G) \mid d(v) \geq 1 \}]. \]

\begin{lem}
Let $k$, $n$, and $r \geq 2$ be positive integers.  If $(G_1,\ldots,G_k)$ is an $\omega$-optimal
$k$-decomposition of $K_n^r$ with $s(G_1)$ maximal, then $P(G_j)$ is complete for $j \geq 2$.
\end{lem}
\begin{proof}
Let $(G_1,\ldots,G_k)$ be am $\omega$-optimal $k$-decomposition of $K_n^r$ with $s(G_1)$ maximal.  Let $j \geq
2$. Take $e \in E(G_j)$.  Then $(G_1 + e,\ldots,G_j - e,\ldots,G_k)$ is a $k$-decomposition of $K_n^r$ with
$s(G_1 + e) > s(G_1)$.  Hence $(G_1 + e,\ldots,G_j - e,\ldots,G_k)$ is not $\omega$-optimal, which implies that
$\omega(G_j - e) < \omega(G_j)$.  Whence every edge of $G_j$ is involved in every maximal clique and thus every
vertex of positive degree is involved in every maximal clique.  Hence $\omega(P(G_j)) = |P(G_j)|$, showing
$P(G_j)$ complete.
\end{proof}

\begin{thm}
Let $k$, $n$, and $r \geq 2$ be positive integers.  Then $\omega(k;K_n^r) \leq n + (r-1){k \choose 2}$ and if $n
\geq (r-1){k \choose 2}$, then $\omega(k;K_n^r) = n + (r-1){k \choose 2}$.
\end{thm}
\begin{proof}
Let $(G_1,\ldots,G_k)$ be a $\omega$-optimal $k$-decomposition of $K_n^r$ with $s(G_1)$ maximal. \newline Set $X
= \displaystyle \bigcup_{j=2}^{k} V(P(G_j))$. By Lemma 12, $P(G_j)$ is complete for $j \geq 2$.  Hence
$P(G_j[X])$ and $P(G_1[X])$ have at most $r-1$ vertices in common for $j \geq 2$.  Thus $|P(G_1[X])| \leq
(r-1)(k-1)$. In particular, $\omega(G_1[X]) = \omega(P(G_1[X])) \leq (r-1)(k-1)$.  We have

\begin{align*}
\omega(k;K_n^r) = \displaystyle \sum_{j=1}^{k} \omega(G_j) &\leq \omega(G_1 - X) + \omega(G_1[X]) + \sum_{j=2}^{k} \omega(G_j)\\
&\leq n - |X| + (r-1)(k-1)+\sum_{j=2}^{k} \omega(G_j) \\
&= n - |X| + (r-1)(k-1)+\sum_{j=2}^{k} |P(G_j)| \\
&\leq n - |X| + (r-1)(k-1)+ |X| + (r-1){k-1 \choose 2} \\
&= n + (r-1){k \choose 2}.
\end{align*}
This proves the upper bound.  To get the lower bound, we generalize a construction in \cite{kostochka}.  The
construction for $n = (r-1){k \choose 2}$ can be extended for each additional vertex by adding all the edges
involving the new vertex to a single hypergraph in the decomposition.  Thus, it will be enough to take care of
the case $n = (r-1){k \choose 2}$.\newline Let $V(K_n^r) = \{ (i,j) \mid 1 \leq i < j \leq k \} \times
\{1,\ldots, r-1 \}$. For $1 \leq t \leq k$, we define a hypergraph $G_t$. Let $V(G_t)$ be the set of vertices of
$K_n^r$ whose names have $t$ in one of the coordinates of the leading ordered pair.  Let $E(G_t)$ be all
$r$-subsets of $V(G_t)$. We have $|V(G_t)| = (r-1)(k-1)$. In addition, the $G_t$ are pairwise edge disjoint
since $i \neq j \Rightarrow V(G_i) \cap V(G_j) \leq r-1$.  Whence $(G_1,\ldots,G_k)$ can be extended to a
$k$-decomposition of $K_n^r$, giving

\[\omega(k;K_n^r) \geq \displaystyle \sum_{j=1}^{k} \omega(G_j) = k(r-1)(k-1) = (r-1){k \choose 2} + (r-1){k \choose
2} = n + (r-1){k \choose 2}.\]

\end{proof}

\end{document}